\documentclass[12pt]{article}
\newtheorem{Th}{Theorem}
\newtheorem{Le}{Lemma}

\newtheorem{St}{Statement}

\begin{document}

\title{Complete intersection theorem and complete nontrivial-intersection theorem for systems of set partitions}

\date{July 10,  2013}
\author{Vladimir Blinovsky}
\date{\small
 Institute for Information Transmission Problems, \\
 B. Karetnyi 19, Moscow, Russia,\\
 Instituto de Matematica e Statistica, USP,\\
Rua do Matao 1010, 05508- 090, Sao Paulo, Brazil\\
vblinovs@yandex.ru}

\maketitle\bigskip

\begin{center}
{\bf Abstract}
\bigskip

We prove the complete intersection theorem and the complete nontrivial-intersection theorem for systems of set partitions. This means that for all positive integers $n$ and $t$ we find the maximum size  of a family of partitions of $n$-element set such that any two partitions from the family have at least $t$ common parts and we also find the maximal size under the additional condition that no $t$ parts appear in all members of the family.\end{center}

\bigskip

{\bf I Introduction}
\bigskip

Let $\Pi (n)$ be the set of partitions of  $[n]$.  Define the intersection of two partitions $p_1 \bigcap p_2 ,\ p_1 ,p_2 \in \Pi (n)$ to be the set of common parts (blocks).  We say that two partitions  $p_1 ,p_2 \in \Pi (n)$
are {\it $t$-intersecting} if the size of their intersection is at least $t$.  A family of partitions is a {\it $t$-intersecting family} if every tow members of it are $t$-intersecting. The collection of $t$-intersecting families of partitions of $[n]$ is denoted by $\Omega (n,t)$. We say that the family of partitions is {\it nontrivially $t$-intersecting family} if it is  $t$-intersecting and fewer than $t$ parts are common to all its members.  The collection of  nontrivially  $t$-intersecting  families of partitions we denote by $\tilde{\Omega} (n,t).$

We say that $i$ is {\it fixed} in a partition  $p\in{\Pi}(n)$ if $\{ i\}$ is a singleton $\{ i\} $ block in $p$. For $p\in\Pi (n)$ let $f(p),\ p\in\Pi (n)$ denote the set of points  fixed by $p.$

Define
\begin{eqnarray*}
&&
M(n,t)=\max\{ |{\cal A}|: {\cal A}\in\Omega (n,t)\} ,\\
&& \tilde{M}(n,t)=\max\{ |{\cal A}|: {\cal A}\in\tilde{\Omega} (n,t)\} .
\end{eqnarray*}
The main result of the present work is obtaining explicit expression for $M(n,t)$ (Theorems~\ref{th33} or Theorem~\ref{th0}) and  $\tilde{M}(n,t)$ (Theorem~\ref{th21}) for all  $n$ and $t$.
The word {\it complete} in the  phrase {\it Complete Intersection Theorem} underline the fact that the problem of determining values $M(n,k)$ and $\tilde{M}(n,k)$ is solved completely for all $n,t$. We also
say that the solution of the above problems are complete.

Let  $B(n)$ be  the number of partitions of the set $[n]$, which is called the {\it Bell number}. Let also $\tilde{B}(n)$
be the number of partitions of the set $[n]$ that do not have singletons.
The Bell numbers satisfy the following relations
\begin{eqnarray}
\label{e1}
&&B(n)=\frac{1}{e}\sum_{i=0}^{\infty}\frac{i^n}{i!} ,\\
&&\tilde{B}(n)=\frac{1}{e}\sum_{i=0}^\infty \frac{i^n}{(i+1)!},\label{e111}
\\ \nonumber &&B(n+1)=\sum_{i=0}^{n}{n\choose i}B(i)
\end{eqnarray}
whereas $\tilde{B}(n)$ satisfies the following relation
\begin{eqnarray}
\label{e2}
&&\tilde{B}(n)=\sum_{i=0}^{n}{n\choose i} (-1)^{n-i} B(i) \label{ee2}
\end{eqnarray}

Define
$$
\gamma (\ell )=\frac{\sum_{i=0}^{n-\ell +1} \tilde{B}\left( n-\frac{\ell +t}{2}+1-i\right){n-\ell +1\choose i}}{\sum_{i=0}^{n-\ell}\tilde{B}\left( n-\frac{\ell +t}{2}-i\right){n-\ell\choose i}} .
$$
Note that, when $\ell$ is fixed,
$$
\gamma (\ell )\to\infty , \hbox{as} \ n\to\infty .
$$
Our first main result is the following theorem.
 \begin{Th}
 \label{th33}

 \begin{eqnarray*}&&
 M (n,t) =\max_{r\in [0,\lfloor (n-t)/2\rfloor ]}|\{ p\in\Pi (n): \biggl| [t+2r]\bigcap f(p)\biggr|\geq t+r\} |.
 \end{eqnarray*}
 \end{Th}\
 It follows from the proof of Theorem~\ref{th33}, that  it can be reformulated as follows:
 \begin{Th}
\label{th0}
Let  $\ell =t+2r$ be the largest number not greater than $n$ satisfying the relation
\begin{equation}
\label{er}
\frac{\ell -t}{2(\ell -1)}\gamma (\ell )\leq 1.
\end{equation}
For this value of $\ell$ we have
\begin{equation}
\label{el}
M(n,t)=\sum_{i=t+r}^{t+2r}{t+2r\choose i}\sum_{j=0}^{n-t-2r}{n-t-2r\choose j}\tilde{B}(n-i-j).
\end{equation}
\end{Th}
Our proof of this theorem is an extension of the ideas from~\cite{3}, where the complete intersection theorem was proved
for a family of  $t$-cycle-intersecting permutations. 

Remark. Each permutation of $[n]$ is determined by the set of cyclic permutations. Cycle-intersection of two permutations is the set of their common cycles. We say that two permutations
are $t$-cycle-intersecting if the size of their intersection is at least $t$.  

It is proved in~\cite{1}  that
$$
M(n,1)=B(n-1) 
$$
and for sufficiently large $n$ in terms of $t$ that 
$$
M(n,t)=B(n-t).
$$
Our theorem completes the solution of the problem of determination of the value $M(n,t)$  for all  $n$ and $t>1.$ .

Let $2^{[n]}$ be the family of subsets of $[n]$ and ${[n]\choose k}$ be the family of $k$-element subsets of $[n].$ We say that a family ${\cal A}\subset 2^{[n]}$ is a $t$-intersecting family if for the arbitrary elements  $a_1 , a_2 \in{\cal A}$ the size of their intersection $\biggl| a_1 \bigcap a_2 \biggr| \geq t .$ Let $I(n,t)$ be the collection of $t$-intersecting 
families  ${\cal A}$ of $[n]$,  $I(n,k,t)$ be the collection of $t$-intersecting $k$-element families  from $[n]$
and $\tilde{I}(n,t),\ \tilde{I}(n,k,t)$ the collection of nontrivially $t$-intersecting families ($ \biggl|\bigcap_{A\in{\cal A}}A\biggr|  <t .$).  Define
$$
\tilde{M}(n,k,t)=\max_{{\cal A}\in\tilde{I}(n,k,t)}|{\cal A}|.
$$
Hilton and Milner proved the next theorem in~\cite{8}.
\begin{Th}
If $n> 2k$, then 
$$
\tilde{M}(n,k,t)={n-1\choose k-1}-{n-k-1\choose k-1}+1.
$$\end{Th}
This theorem was proved by Frankl~\cite{9} for $t>1$.
\begin{Th}
There exists $n_0 (n,k)$ such that if $n>n_0 (n,k)$, then
\begin{itemize}
\item If $t+1\leq k\leq 2t+1$, then $\tilde{M}(n,k,t)=|\nu_1 (n,k,t)|$, where
$$
\nu_1 (n,k,t)=\left\{ V\in{[n]\choose k}:\ \biggl| [t+2]\bigcap V\biggr|\geq t+1\right\} ,
$$
\item If $k>2t+1$, then $\tilde{M}(n,k,t)=|\nu_2 (n,k,t)|,$ where
\begin{eqnarray*}
\nu_2 (n,k,t)&=&\left\{ v\in{[n]\choose k}: [t]\subset V,\ V\bigcap [t+1,k+1]\neq\emptyset\right\}\\
&\bigcup& \left\{ [k+1]\setminus\{ i\}:\ i\in [t]\right\} .
\end{eqnarray*}
\end{itemize}
\end{Th}
In~\cite{7}, the problem of determining $\tilde{M}(n,k,t)$ was solved completely for all $n,k,t$:
\begin{Th}\mbox{}\\
\label{th66}
\begin{itemize}
\item If $2k-t<n\leq (t+1)(k-t+1)$, then
$$
\tilde{M}(n,k,t)=M(n,k,t);
$$
\item If $(t+1)(k-t+1)<n$ and $k\leq 2t+1$, then
$$
\tilde{M}(n,k,t)=|\nu_1 (n,k,t)|;
$$
\item If $ (t+1)(k-t+1)<n$ and $k>2t+1$, then
$$
\tilde{M}=\max\{ |\nu_1 (n,k,t)|,|\nu_2 (n,k,t)|\}.
$$
\end{itemize}
\end{Th}
Note also that the value $M(n,k,t)$ was determined for all $n,k,t$ by Ahlswede and Khachatrian in the paper~\cite{11}. 

Before formulating our second main result, let's make some additional definitions.

\begin{eqnarray*}
&& {\cal H}_i =\biggl\{ H\in{[t+i]\choose t+1}: [t]\subset H\biggr\} \\
&&\bigcup \biggl\{ H\in{[t+i]\choose t+i-1}: [t+1,t+i]\subset H\biggr\} .
\end{eqnarray*}
For ${\cal C}\subset 2^{[n]},$ denote by $W({\cal C})$ the minimal upset containing ${\cal C}$ and by $M({\cal C})$ the set of its minimal elements. Denote by $U({\cal C})$ the set of partitions 
that has $W({\cal C})$ as the family of sets of fixed elements.

 Our second main result of this work is  the following Theorem which completely determines
 $\tilde{M}(n,t)$ for all $n,t.$ 
 \begin{Th}\mbox{}\\
 \label{th21}
 \begin{itemize}
 \item If 
 $$
 \max\left\{ \ell=t+2r: \frac{\ell-t}{2(\ell -1)}\gamma (\ell )\leq 1\right\} >t,
 $$
 then
 $$
 \tilde{M}(n,t)=M(n,t);
 $$
 \item
 If
 $$
 \max\left\{ \ell=t+2r: \frac{\ell-t}{2(\ell -1)}\gamma (\ell )\leq 1\right\} =t ,
 $$
 then
 $$
 \tilde{M}(n,t)=\max\{\nu_1 (n,t) ,\nu_2(n,t)\} ,
 $$
 where
 $$
 \nu_i (n,t)=\sum_{S\in W({\cal H}_i )}\tilde{B}(n-|S|) .
 $$
 \end{itemize}
 \end{Th}
 \bigskip
 
{\bf II Proof of Theorem~\ref{th33}}
\bigskip

Define the {\it fixing} procedure $F(i,j,p)$ for $ i\neq j$ over the set of partitions $p\in{\cal P}(n)$:
$$
F(i,j,p)=\left\{\begin{array}{ll}
(p\setminus p_i )\bigcup \{ \{ i\} ,p_i \setminus\{i\}\},& j\in p_i\\
p , \hbox{otherwise}
\end{array},
\right.
$$
where $p_i$ is the part of $p$ that contains $i.$

The fixing operator on the family ${\cal A}\subset\Omega (n,t)$
is defined as follows $(p\in{\cal A})$
$$
F(i,j,p,{\cal A})=\left\{\begin{array}{ll}
F(i,j,p),\ & F(i,j,p )\not\in {\cal A},\\
p, & F(i,j,p) \in{\cal A} .
\end{array}\right.
$$
Finally define the operator
$$
{\cal F} (i,j,{\cal A})=\{ F(i,j,p,{\cal A}); p\in{\cal A}\}.
$$
It is easy to see that the fixing operator ${\cal F}(i,j,{\cal A})$ preserves the size of ${\cal A}$ and its $t$-intersecting property. At last note that making shifting operations a finite number of times for different values of $i$ and $j$ allows us to obtain the compressed set {\cal A} with the following property: for all $i\neq j\in [n]$,
$$
{\cal F}(i,j,{\cal A}) ={\cal A}.
$$
It also has the property, that an arbitrary pair of partitions $p_1 ,p_2$ from the compressed set ${\cal A}$ intersected by at least $t$
fixed points.

Next define the usual  {\it shifting} procedure $L (v,w,p)$ for $ 1\leq v<w\leq n$ as follows. 
Let $p=\{ \{ j_1 ,\ldots ,j_{q-1},v,j_{q+1},\ldots ,j_s \} ,\ldots , \{ w\} ,\pi_1 ,\ldots ,\pi_c \} \in{\cal A}$, then
$$
L(v,w,p)=\{ \{j_1 ,\ldots ,j_{q-1},w,j_{q+1},\ldots ,j_s\} ,\ldots ,\{ v\} ,\pi_1 ,\ldots ,p_s \} .
$$
If $p\in{\cal A}$ does not fix $w$, then we set
$$
L(v,w, p)=p .
$$
Now define the shifting operator $L(v,w, p, {\cal A})$ as follows
$$
L(v,w,p,{\cal A})=\left\{\begin{array}{ll}
L(v,w,p),\ & L (v,w,p)\not\in{\cal A},\\
p, & L(v,w,p)\in{\cal A} .
\end{array}
\right.
$$
At last define the operator ${\cal L}(v,w,{\cal A}) :$
$$
{\cal L}(v,w,{\cal A})=\{ L(v,w, {\cal A});\ p\in{\cal A}\} .
$$
It is easy to see that the operator ${\cal L}(v,w,{\cal A})$ does not change the size of ${\cal A}$ and it preserve the $t$-intersecting property. Later we will show, proving the Statement~\ref{s1}, that this operator also preserves the  nontrivially $t$-intersecting property. Also it is easy to see that after a finite number of operations we come  to the compressed  $t$-intersecting set ${\cal A}$ of the size $\tilde{M}(n,t)$ for which
$$
L(v,w,{\cal A})={\cal A}\  \hbox{for}\ 1\leq v<w\leq n 
$$
and to the property that each pair of partitions of ${\cal A}$ is $t$-intersected by fixed elements.
Next we consider only such sets ${\cal A}$. 
We denote the collection of fixed compressed  $t$- intersecting families of partitions by $L\Omega (n,t)$ and the collection of fixed compressed nontrivially $t$-intersecting families of partitions by $ (L\tilde{\Omega}(n,t))$.
Note that such family ${\cal A}$
have the property that all partitions of ${\cal A}$ have $s$ common parts if and only if  $\biggl|\bigcap_{p\in{\cal A}}f(p)\biggr| =s.$
We assume that all families of partitions considered next are left compressed.

Let ${\cal D}(v,w,{\cal A})$ be the same operator as ${\cal L}(v,w,{\cal A})$ but only with the condition $v\neq w$.

We need the following
\begin{Le}
\label{le2}
If $|{\cal A}|=M(n,t),$ 
$$
D(v,w,{\cal A})={\cal A},\ \hbox{ for all}\ v, w \in [\ell ] 
$$
and for $\ell =t+2r$ 
\begin{equation}
\label{e99}
\frac{\ell +1}{\frac{\ell -t}{2}+1}>\gamma (\ell +2),
\end{equation}
then  ${\cal D}(v,w,{\cal A})={\cal A},\  \hbox{for all}\  v,w\in [\ell +2] .$
\end{Le}
 Suppose  that
 $$
 |{\cal A}| =M(n,t)$$ and ${\cal A}$ is invariant under shifting and fixing operators . Assume also that 
 $$
 {\cal D}(v,w,{\cal A})={\cal A} \ \hbox{ for all}\ v,w\in [\ell ],
 $$
 but
 $$
 {\cal A}\neq {\cal D}(v,\ell +1,{\cal A})
 $$
 for some $v\in [\ell ].$
 
 We set
 $$
 {\cal A}^\prime =\{ p\in{\cal A}: D(v,\ell +1 ,p)\not\in{\cal A}  , v\in [\ell]\} 
 $$  
We identify the set of binary $n$-tuples with the family of subsets of $[n].$

Define
$$
B({\cal A})=\{ f(p);\ p\in{\cal A}\}\subset 2^{[n]} .
$$
It is easy to see that the set $B({\cal A} )$ is an upper ideal under the inclusion order. Denote by $M({\cal A})$ the set of minimal elements of $B({\cal A})$.

Let
$$
s^+ (M({\cal A}))=\max_{M\in M({\cal A})}s^+ (M) ,
$$
where
$$
s^+ (M)=\max_{i\in M}i .
$$
We also need one more lemma.
\begin{Le}
\label{le3}
If  $|{\cal A}|=M(n,t)$ and $ s^+ (M({\cal A}))=\ell$, then
\begin{equation}
\label{e89}
\frac{\ell -t}{2(\ell -1)}\gamma (\ell )\leq 1.
\end{equation}
\end{Le}
Later we will show that there exists a unique $\ell$ that satisfies inequalities~(\ref{e99}) and (\ref{e89}).
From this follows the statement of Theorem \ref{th33}.

First we will prove~\ref{le2}.
Assume that ${\cal A}^\prime \neq\emptyset .$

Let
$$
{\cal A}(i)=\left\{ p\in {\cal A}^\prime :\ \biggr| f(p)\bigcap [\ell ]\biggr| =i\right\} .
$$
It follows that ${\cal A}(i)\neq\emptyset $ for some $i\in [\ell ]$. 

Let
$$
{\cal A}^\prime (i)= \left\{ f(p)\bigcap [\ell +2 ,n]:\ p\in {\cal A}(i)\right\} .
$$
Define
$${\cal B}(i)=\left\{ p\in\Pi (n): \biggl| f(p)\bigcap [\ell ]\biggr| =i-1,\ \ell +1\in p,\ f(p)\bigcap [\ell +2,n]\in{\cal A}^\prime (i )\right\} .
$$
We have
$$
|{\cal A} (i)|={\ell\choose i}\sum_{m\in{\cal A}^\prime (i)}\tilde{B}(n-i- 1-|m|),
$$
$$
|{\cal B}(i)|={\ell -1\choose i-1}\sum_{m\in{\cal A}^\prime (i)}\tilde{B} (n-i-1-|m|) 
$$
and, for any $i\in [\ell],$ we have
$$
{\cal C}(i)=({\cal A}\setminus {\cal A}(i))\bigcup{\cal B} (\ell +t-i)\in\Omega (n,t).$$

Next we will demonstrate that if ${\cal A}(i)\neq\emptyset $ and  $i\neq \frac{\ell +t}{2} $, then
\begin{equation}
\label{e4}
\max\{ |{\cal C}(i)|, |{\cal C}(\ell +t-i)|\} >|{\cal A}|
\end{equation}
which contradicts the maximality of ${\cal A}$.

If~(\ref{e4}) is not valid, then
\begin{eqnarray*}
&& {\ell\choose i-1}\sum_{m\in{\cal A}^\prime (i)}\tilde{B} (n-i-1-|m|)\\
&&\leq{\ell\choose\ell +t-i}\sum_{r\in{\cal A}^\prime (\ell +t-i)}\tilde{B}(n-(\ell +t-i)-1-|m|) ,\\
&&
{\ell\choose\ell +t-i-1}\sum_{m\in{\cal A}^\prime (\ell +t -i)}\tilde{B} (n-(\ell +t-i)-1-|m|)\\
&&\leq {\ell\choose i}\sum_{m\in{\cal A}^\prime (i)}\tilde{B} (n-i-1-|m|) .
\end{eqnarray*}
We have ${\cal A}(i)\neq\emptyset$ hence ${\cal A} (\ell +t-i )\neq\emptyset  $
and
$$
i(\ell +t -i)\leq (\ell -i +1) (i+1-t ).
$$
Since $t\geq 2$, the last inequality is false. This contradiction shows that ${\cal A}(i)=\emptyset$ for $i\neq\frac{\ell +t}{2} .$

Now suppose $2|(\ell +t)$. We will demonstrate that if~(\ref{e99}) is true, then
${\cal A}\left(\frac{\ell +t}{2}\right) =\emptyset .$

We have 
$$
\biggl|{\cal A}\left(\frac{\ell +t}{2}\right)\biggr| ={\ell\choose \frac{\ell +t}{2}}\sum_{m\in{\cal A}^\prime\left(\frac{\ell +t}{2}\right)}\tilde{B} \left(n-\frac{\ell+t}{2}-1 -|m|\right).
$$
Now we will introduce the family ${\cal C}\subset\Pi (n)$. Its elements are permutations $p$ which satisfy the following conditions
\begin{eqnarray*}
&& \ \biggl| f(p)\bigcap [\ell ]\biggr| =\frac{\ell +t}{2}-1 ,\\
&&\ \{\ell +1 ,n\}\subset f(p) ,\  f(p)\bigcap [\ell +2 ,n]\in{\cal A}^\prime \left(\frac{\ell +t}{2}\right) .
\end{eqnarray*}

Define 
$$
{\cal G}=\left\{ \left({\cal A} \setminus \left\{ p\in{\cal A}\left(\frac{\ell +t}{2}\right) :\ \{n\}\not\in p\right\} \right)\bigcup {\cal C}\right\} .
$$
It is easy to see that$$
{\cal G}\subset\Omega (n,t).
$$
Next we will demonstrate that if ${\cal A}\left( \frac{\ell +t}{2}\right)\neq \emptyset$ and~(\ref{e99}) is true, then
the maximality of ${\cal A}$ is contradicted because
\begin{equation}
\label{e10}
|{\cal G}|>|{\cal A}| .
\end{equation}
We have
$$
|{\cal C}|={\ell\choose\frac{\ell +t}{2}-1}\sum_{m\in{\cal A}^\prime \left(\frac{\ell +t}{2}\right)
 ,\{ n\}\in m}
\tilde{B}\left( n-\frac{\ell +t}{2}-1-|m|\right) .
$$
Inequality~(\ref{e10}) is equivalent to
\begin{eqnarray*}
&&{\ell\choose\frac{\ell +t}{2}-1}\sum_{m\in{\cal A}^\prime \left(\frac{\ell +t}{2}\right),\ \{ n\} \in m}\tilde{B}\left(n-\frac{\ell +t}{2}- 1-|m|\right) \\
&>& {\ell\choose\frac{\ell +t}{2}}\sum_{m\in{\cal A}^\prime \left(\frac{\ell +t}{2}\right) ,\ \{n\}\not\in m}
\tilde{B}\left(n-\frac{\ell +t}{2}-1-|m|\right) \\
&=&{\ell\choose\frac{\ell +t}{2}}\sum_{m\in{\cal A}^\prime \left(\frac{\ell +t}{2}\right)} \tilde{B}\left( n-\frac{\ell +t}{2}-1-|m|\right) \\
&-& \sum_{m\in{\cal A}^\prime \left(\frac{\ell +t}{2}\right) ,\ \{n\}\in m}\tilde{B}\left( n-\frac{\ell +t}{2}-1-|m|\right) .
\end{eqnarray*}
From here we have
\begin{eqnarray*}
&&
{\ell+1\choose\frac{\ell +t}{2}}\sum_{m\in{\cal A}^\prime \left(\frac{\ell +t}{2}\right) ,\ \{ n\}\in m}\tilde{B} \left( n-\frac{\ell +t}{2} -1-|m|\right) \\
&>& {\ell\choose\frac{\ell +t}{2}}\sum_{m\in{\cal A}^{\prime} \left(\frac{\ell +t}{2}\right)}\tilde{B} \left( n-\frac{\ell +t}{2} -1-|m|\right) .
\end{eqnarray*}
Hence
$$
\frac{\ell +1}{\frac{\ell -t}{2}+1}> \beta_1 (\ell)\stackrel{\Delta}{=}\frac{\sum_{m\in{\cal A}^\prime \left(\frac{\ell +t}{2}\right)} \tilde{B} \left( n-\frac{\ell +t}{2}-1-|m|\right)}{\sum_{m\in{\cal A}^\prime \left (\frac{\ell +t}{2}\right),\ \{ n\} \in m}\tilde{B} \left( n-\frac{\ell +t}{2} -1-|m|\right)} .
$$
 Let's prove that
 \begin{equation}
 \label{e0}
 \gamma (\ell +2)\geq\beta_1 (\ell ) .
 \end{equation}
From here it follows that~(\ref{e99}) is true. Taking into account the condition from Lemma~\ref{le2} that $|{\cal A}|$ is maximal we come to contradiction of this maximality.  Thus to complete the proof of Lemma~\ref{le2} we need to prove equality~(\ref{e0}).

 In order to show this, we state the validity of the following inequality 
\begin{eqnarray*}
&&\gamma (\ell +2)=\frac{\sum_{i=0}^{n-\ell -1}\tilde{B} \left(n-\frac{\ell+t}{2}-i\right){n-\ell-1\choose i}}{\sum_{i=0}^{n-\ell -2}\tilde{B} \left(n-\frac{\ell+t}{2}-i-1\right){n-\ell-2\choose i}}\geq\frac{\sum_{m\in{\cal A}^\prime \left(\frac{\ell +t}{2}\right)} \tilde{B} \left( n-\frac{\ell +t}{2}-|m|\right)}{\sum_{m\in{\cal A}^\prime \left (\frac{\ell +t}{2}\right),\ \{ n\} \in m}\tilde{B} \left( n-\frac{\ell +t}{2} -|m|\right)}=\beta_1 \nonumber\end{eqnarray*}
or
\begin{eqnarray}
\label{s111}&&\\
&&\frac{\sum_{i=0}^{n-\ell -2}\tilde{B} \left(n-\frac{\ell+t}{2}-i\right){n-\ell-2\choose i}}{\sum_{i=0}^{n-\ell -2}\tilde{B} \left(n-\frac{\ell+t}{2}-i-1\right){n-\ell-2\choose i}}\geq\frac{\sum_{m\in{\cal A}^\prime \left(\frac{\ell +t}{2}\right),\ \{n\}\not\in m} \tilde{B} \left( n-\frac{\ell +t}{2}-|m|\right)}{\sum_{m\in{\cal A}^\prime \left (\frac{\ell +t}{2}\right),\ \{ n\} \in m}\tilde{B} \left( n-\frac{\ell +t}{2} -|m|\right)}. \nonumber\end{eqnarray}
Validity of this inequality follows from the next consideration. Family ${\cal A}$ is compressed under fixing operator. Set of minimal elements $M({\cal A})\subset 2^{[\ell]}$ in the set $B({\cal A})$ has $t$- intersection property. 

Hence  $\{f(p)\bigcap [\ell+2,n];\ p\in {\cal A}\} =2^{[\ell+2,n]}$ and inequality~(\ref{s111}) can be written as follows:
\begin{eqnarray*}&&
\frac{\sum_{i=0}^{n-\ell -2}\tilde{B} \left(n-\frac{\ell+t}{2}-i\right){n-\ell-2\choose i}}{\sum_{i=0}^{n-\ell -2}\tilde{B} \left(n-\frac{\ell+t}{2}-i-1\right){n-\ell-2\choose i}}\geq \frac{\sum_{i=0}^{n-\ell -2}\tilde{B} \left(n-\frac{\ell+t}{2}-i-1\right){n-\ell-2\choose i}}{\sum_{i=0}^{n-\ell -2}\tilde{B} \left(n-\frac{\ell+t}{2}-i-2\right){n-\ell-2\choose i}}.
\end{eqnarray*}

 or
 \begin{eqnarray}&&\label{f4}
\frac{\sum_{i=0}^{n-\ell -2}\tilde{B} \left(n-\frac{\ell+t}{2}-i-2\right){n-\ell-2\choose i}}{\sum_{i=0}^{n-\ell -2}\tilde{B} \left(n-\frac{\ell+t}{2}-i-1\right){n-\ell-2\choose i}}\geq \frac{\sum_{i=0}^{n-\ell -2}\tilde{B} \left(n-\frac{\ell+t}{2}-i-1\right){n-\ell-2\choose i}}{\sum_{i=0}^{n-\ell -2}\tilde{B} \left(n-\frac{\ell+t}{2}-i\right){n-\ell-2\choose i}}
\end{eqnarray}
or, using identity
$
\tilde{B}(m)=B(m)-\tilde{B}(m+1),$ we obtain from the inequality~(\ref{f4}) inequality
\begin{eqnarray}&& \label{kl}
\frac{\sum_{i=0}^{n-\ell -2}B \left(n-\frac{\ell+t}{2}-i-2\right){n-\ell-2\choose i}}{\sum_{i=0}^{n-\ell -2}\tilde{B} \left(n-\frac{\ell+t}{2}-i-1\right){n-\ell-2\choose i}}\geq \frac{\sum_{i=0}^{n-\ell -2}B \left(n-\frac{\ell+t}{2}-i-1\right){n-\ell-2\choose i}}{\sum_{i=0}^{n-\ell -2}\tilde{B} \left(n-\frac{\ell+t}{2}-i\right){n-\ell-2\choose i}}
\end{eqnarray}
 Validity of last inequality follows from Holley's correlation inequality. Let's $\Delta\subseteq 2^{[n]}$ be finite distributive lattice and measures   $\mu_1,\ \mu_2 :\ \Delta\to  R_+\bigcup\{0\},\ \sum_{A\in\Delta}\mu_1 (A)=\sum_{A\in\Delta}\mu_2(A)$ satisfy  FKG condition  
 \begin{equation}
 \label{l4}
 \mu_1 (A)\mu_2 (B)\leq\mu_1\left(A\bigcup B\right)\mu_2 \left(A\bigcap B\right),\ A,B\in \Delta.
 \end{equation}
 Then for an arbitrary nondecreasing nonnegative function $f:\ \Delta\to R_+\bigcup\{0\},\  A\subset B\to f(A)\geq f(B),$
 the Holley's inequality stand
 $$
  \sum_{A\in\Delta}\mu_1(A)f(A) \geq \sum_{A\in \Delta}\mu_2 (A)f(A).
  $$
   We choose 
   \begin{eqnarray*}
   &&
   \Omega =\{ [k],\ k\in [n-\ell-2]\},\ f([k])=\frac{B\left(n-\frac{\ell+t}{2}-k-2\right)}{\tilde{B}\left(n-\frac{\ell+t}{2}-k-1\right)},\\
   &&
   \mu_{1}([k])=\frac{\tilde{B}\left(n-\frac{\ell+t}{2}-k-1\right){n-\ell-2\choose k}}{\sum_{i=0}^{n-\ell -2}\tilde{B} \left(n-\frac{\ell+t}{2}-i-1\right){n-\ell-2\choose i}},\\
   &&
   \mu_{2}([k])=\frac{\tilde{B}\left(n-\frac{\ell+t}{2}-k\right){n-\ell-2\choose k}}{\sum_{i=0}^{n-\ell -2}\tilde{B} \left(n-\frac{\ell+t}{2}-i\right){n-\ell-2\choose i}}.
   \end{eqnarray*}
   FKG condition for $\mu_1,\ \mu_2$ follows from the validity if FKG condition for Bell numbers $B(k)$ proved in~\cite{14}.
   
   We need to check the monotonicity of $f([k]).$ If $f([k+1])\geq f([k]),$ then
   $$
   \frac{B\left(n-\frac{\ell+t}{2}-k-2\right)}{\tilde{B}\left(n-\frac{\ell+t}{2}-k-1\right)} \geq \frac{B\left(n-\frac{\ell+t}{2}-k-1\right)}{\tilde{B}\left(n-\frac{\ell+t}{2}-k\right)}
   $$
   or 
    identity $B(m)=\tilde{B}(m)+\tilde{B}(m+1),$ we  reduce last inequality to the so -called log - convexity condition of $\tilde{B}(m)$:
   $$
   \tilde{B}^2\left(n-\frac{\ell+t}{2}-k-1\right) \leq    \tilde{B}\left(n-\frac{\ell+t}{2}-k-2\right)     \tilde{B}\left(n-\frac{\ell+t}{2}-k\right).
   $$
   In Appendix we prove that  $\tilde{B}\left( n-\frac{\ell+t}{2}-k-1\right)$  satisfies this inequality when $n-\frac{\ell +t}{2}-k-1>3$. This inequality is true for all $k\in[n-\ell-2]$ and  $r>2.$
   
Assume  at first that $r>2.$ Then, using Holley's inequality we have
   \begin{eqnarray*}&&
\sum_{k=0}^{n-\ell-2}\mu_1 ([k])f([k+1])=\frac{\sum_{k=0}^{n-\ell -2}B \left(n-\frac{\ell+t}{2}-k-2\right){n-\ell-2\choose i}}{\sum_{i=0}^{n-\ell -2}\tilde{B} \left(n-\frac{\ell+t}{2}-i-1\right){n-\ell-2\choose i}}\\&& \geq \sum_{k=0}^{n-\ell-2}\mu_2 ([k])f([k+1]) \geq \sum_{k=0}^{n-\ell-2}\mu_2 ([k])f([k])=\frac{\sum_{i=0}^{n-\ell -2}B \left(n-\frac{\ell+t}{2}-i-1\right){n-\ell-2\choose i}}{\sum_{i=0}^{n-\ell -2}\tilde{B} \left(n-\frac{\ell+t}{2}-i\right){n-\ell-2\choose i}}
\end{eqnarray*}
This proves inequality~(\ref{f4}) when $r> 2$. When $r\leq 3$, we use identities~\cite{20} (both identities are actually the same, we just write them in form in which they we will use them later):
\begin{eqnarray*}
&&
{n-t-2r-2\choose i}=\sum_{k=0} (-1)^k {q\choose k}{n-t-2r-2 +q-k\choose i+q},\ q\geq 0,\\&&
{n-t-2r-2\choose a+i-r-2}=\sum_{k=0}(-1)^k {r+2-a\choose k}{n-t-r-a-k\choose i}.
\end{eqnarray*}
We have $(a\leq 2+r):$
\begin{eqnarray*}&&
\sum_{i=0}^{n-t-2r -2}\tilde{B} \left(n-t-r-i-a\right){n-t-2r-2\choose i} = \sum_{i=0}^{n-t-2r -2}\tilde{B} \left(r+2-a+i\right){n-t-2r-2\choose i}\\
&&=  \sum_{k=0} (-1)^k {r+2-a\choose k}\sum_{i=0}\tilde{B} \left(r+2-a+i\right) {n-t-r-a-k\choose  i+r+2-a}\\
&&=   \sum_{k=0} (-1)^k {r+2-a\choose k}\sum_{i=r+2-a}^{n-t-r-a-k}\tilde{B} (i) {n-t-r-a-k\choose  i}\\&&
=  \sum_{k=0} (-1)^k {r+2-a\choose k}\left( B(n-t-r-a-k)-\sum_{i=0}^{r+1-a}\tilde{B} (i) {n-t-r-a-k\choose  i}\right)\\
&&=\sum_{k=0} (-1)^k {r+2-a\choose k}B(n-t-r-a-k) -\sum_{i=0}^{r+1-a}\tilde{B}(i)\sum_{k=0} (-1)^k {r+2-a\choose k}{n-t-r-a-k\choose  i}\\
&&=\sum_{k=0}^{r+2-a} (-1)^k {r+2-a\choose k}B(n-t-r-a-k).
\end{eqnarray*}
Using last identity we can rewrite inequality~(\ref{f4}) as follows
\begin{eqnarray*}
&&
\frac{\sum_{k=0}^{r+1} (-1)^k {r+1\choose k}B(n-t-r-1-k)}{\sum_{k=0}^{r} (-1)^k {r\choose k}B(n-t-r-2-k)}\leq \frac{\sum_{k=0}^{r+2} (-1)^k {r+2\choose k}B(n-t-r-k)}{\sum_{k=0}^{r+1} (-1)^k {r+1\choose k}B(n-t-r-1-k)}
\end{eqnarray*}
Denoting $S(r)=\sum_{k=0}^{r} (-1)^k {r\choose k}B(n-t-r-2-k)$, we need to show the validity of inequality
$$
S^2(r+1)\leq S(r+2)S(r).
$$
We write this inequality for $r=0,1,2,3$ making some transformations:
\begin{eqnarray*}
&&
\frac{B(n-t-1)}{B(n-t-2)}\leq \frac{B(n-t)}{B(n-t-1)},\ r=0;\\
&&
\frac{B(n-t-2)-B(n-t-3)}{B(n-t-3)-B(n-t-4)}\leq\frac{B(n-t-1)-B(n-t-2)}{B(n-t-2)-B(n-t-3)},\ r=1;\ \\
&&
\frac{B(n-t-3)-2B(n-t-4)+B(n-t-5)}{B(n-t-4)-2B(n-t-5)+B(n-t-6)}\\
&&\leq \frac{B(n-t-2)-3B(n-t-3)+3B(n-t-4)-B(n-t-5)}{B(n-t-3)-3B(n-t-4)+3B(n-t-5)-B(n-t-6)},\ r=2;\\
&&\frac{B(n-t-4)-3N(n-t-5)+3B(n-t-6)-B(n-t-7)}{B(n-t-5)-3B(n-t-6)+3B(n-t-7)-B(n-t-8)}\\
&&\leq \frac{B(n-t-3)-4B(n-t-4)+6B(n-t-5)+4B(n-t-6)-B(n-t-7)}{B(n-t-4)-4B(n-t-5)+6B(n-t-6)-4B(n-t-7)+B(n-t-8)},\ r=3.
\end{eqnarray*}
  Next we will prove Lemma~\ref{le3}. 
  
  Define
$$
M_0 ({\cal A})=\{ E\in M({\cal A}); s^+ (E)=s^+ (M({\cal A}))=\ell \}
$$
and
$$
M_1 ({\cal A})=M({\cal A})\setminus M_0 ({\cal A}) .
$$
It is easy to see that, for
$E_1 \in M_0 ({\cal A})$ and $E_2 \in M_1 ({\cal A}),$
$$
\biggl| (E_1 \setminus\{\ell\} )\bigcap E_2 \biggr|\geq t
$$
and for $E_1 ,E_2 \in M_0 ({\cal A})$ and $\ \biggl| E_1 \bigcap E_2 \biggr| =t,$ 
$$
|E_1 |+ |E_2 |=\ell +t. 
$$
Set
$$
M_0 ({\cal A})=\bigcup_i R(i) ,
$$
where
$$
R(i)=M_0 ({\cal A})\bigcap {[n]\choose i} .
$$
Define
$$
R^\prime (i)=\{ E\setminus\{\ell\} ;\ E\in R(i)\} .
$$
Next we are going to prove that if~(\ref{e99}) is not true, then $R(i)=\emptyset .$

Suppose that $R(i)\neq\emptyset$ for some $i.$ At first, assume that $i\neq\frac{\ell +t}{2} .$

Define
\begin{eqnarray*}
&&F_1 =M_1 ({\cal A})\bigcup (M_0 ({\cal A})\setminus \left( R(i)\bigcup R(\ell +t-i)\right) )\bigcup R^\prime (i),\\
&& F_2 =M_1 ({\cal A})\bigcup (M_0 ({\cal A})\setminus \left( R(i)\bigcup R(\ell +t-i)\right) )\bigcup R^\prime (\ell +t-i) .
\end{eqnarray*}
It is easy to see that for $E_1 ,E_2 \in F_i $ we have $\ \biggl| E_1 \bigcap E_2 \biggr|\geq t$ and thus $U(F_1 ),\ U(F_2 )\in\Omega (n,t) .$
We are going to show that if $R(i)\neq\emptyset ,$ then
\begin{equation}
\label{e90}
\max\{ |U(F_1 )|, |U(F_2 )|\} >|{\cal A}| 
\end{equation}
which gives us a contradiction.

We have
\begin{equation}
\label{e11}
|{\cal A}\setminus U(F_1 )|=|R(\ell +t-i)|\sum_{j=0}^{n-\ell}{n-\ell\choose j}\tilde{B} (n-\ell -t+i-j)
\end{equation}
and
\begin{equation}
\label{e112}
|U(F_1 )\setminus{\cal A}|=|R(i)|\sum_{j=0}^{n-t}{n-\ell\choose j}\tilde{B} (n-i-j+1) .
\end{equation}
Also
\begin{eqnarray}
\label{e13}
&&|{\cal A}\setminus U(F_2 )|=|R(i)|\sum_{j=0}^{n-\ell}{n-\ell\choose j}\tilde{B} (n-i-j) ,\\
\label{e14}
&&
|U(F_2 )\setminus{\cal A}|=|R(\ell +t-i)|\sum_{j=0}^{n-\ell}\tilde{B}(n-\ell -t+i-j+1) .
\end{eqnarray}
If~(\ref{e90}) is not true, then from~(\ref{e11})-(\ref{e14}) it follows that
\begin{eqnarray*}
&&
|R(i)|\sum_{j=0}^{n-\ell}{n-\ell\choose j}\tilde{B}(n-i-j+1)\\
&&\leq |R(\ell +t-i)|\sum_{j=0}^{n-\ell}{n-\ell\choose j}\tilde{B} (n-\ell -t+i-j) ,\\
&&
|R(\ell +t-i)|\sum_{j=0}^{n-\ell}{n-\ell\choose j}\tilde{B} (n-\ell -t +i-j+1)\\
&&\leq |R(i)|\sum_{j=0}^{n-\ell}{n-\ell\choose j}\tilde{B} (n-i-j ).
\end{eqnarray*}
These inequalities couldn't be valid together due to monotonicity of $\tilde{B}(n).$

Now consider the case $i=\frac{\ell +t}{2} .$ We are going to prove that if inequality~(\ref{e89}) is not true, then
$R\left(\frac{\ell +t}{2}\right) =\emptyset .$ Simple averaging argument shows that there exists $i\in [\ell -1]$ and $Z\subset R^\prime \left(\frac{\ell +t}{2}\right)$  such that $i\in E$ for all $E\in Z$ and
\begin{equation}
\label{e777}
|Z|\geq\frac{\ell -t}{2(\ell -1)}\biggl| R^\prime \left(\frac{\ell +t}{2}\right) \biggr| .
\end{equation}
Because
$\biggl| E_1 \bigcap E_2 \biggr|\geq t$ when $ E_1 ,E_2 \in Z$ and $R(i)=\emptyset $ when $ i\neq\frac{\ell +t}{2} $ we have  for all $E_1 ,E_2 \in D$, where
$$
D=\left( M({\cal A})\setminus R\left(\frac{\ell +t}{2}\right)\right)\bigcup Z,
$$
we have $|E_1 \bigcap E_2 |\geq t.$ Hence $W(D)\in\Omega (n,t)$ and now we have to show that, if~(\ref{e89}) is not true, then
\begin{equation}
\label{e23}
|W(D)|>|{\cal A}| .
\end{equation}
Consider the partition
\begin{eqnarray*}
&&{\cal A}=W(M({\cal A}))=S_1 \bigcup S_2 ,\\
&& S_1 =W\left( M({\cal A})\setminus R\left(\frac{\ell +t}{2}\right)\right) ,\\
&& S_2 =W\left( R\left(\frac{\ell +t}{2}\right)\right)\setminus S_1
\end{eqnarray*}
and the partition
\begin{eqnarray*}
&&W(D)=S_1 \bigcup S_3 ,\\
&& S_3 =W(D)\setminus S_1 .
\end{eqnarray*}
One can see that~(\ref{e23}) is equivalent to
$$
|S_3 |> |S_2 | .
$$
It is easy to show that
\begin{eqnarray}
\label{e778}
&& |S_2 |=\biggl| R\left(\frac{\ell +t}{2}\right) \biggr| \sum_{j=0}^{n-\ell}{n-\ell\choose j}\tilde{B} \left( n-\frac{\ell +t}{2}-j\right) ,\\
&& |S_3 |=|Z|\sum_{j=0}^{n-\ell +1}{n-\ell +1\choose j}\tilde{B} \left( n-\frac{\ell +t}{2}+1-j\right) .
\nonumber
\end{eqnarray}
Using~(\ref{e777}) and~(\ref{e778}) we conclude  that
\begin{eqnarray*}
&&\frac{\ell -t}{2(\ell -1)}\biggl| R\left(\frac{\ell +t}{2}\right)\biggr|\sum_{j=0}^{n-\ell +1}{n-\ell +1\choose j}\tilde{B}\left( n-\frac{\ell +t}{2}+1-j\right)\\
&>& \biggl| R\left(\frac{\ell +t}{2}\right)\biggr| \sum_{j=0}^{n-\ell}{n-\ell\choose j}\tilde{B}\left( n-\frac{\ell +t}{2}-j\right) .
\end{eqnarray*}
But from here follows the contradiction of the maximality of ${\cal A}$. 

Thus~(\ref{e89}) holds.

Now we rewrite inequality~(\ref{e99}) as follows
$$
\ell +2< t+2\frac{t-1}{\gamma (\ell +2)-2} ,
$$
and inequality~(\ref{e89}) as
$$
\ell \leq t+2\frac{t-1}{\gamma (\ell )-2} .
$$
It is left for us to show that the function
$$
\varphi (\ell )= t-\ell+2\frac{t-1}{\gamma (\ell )-2}
$$
does not change its sign in the interval $[t,n]$ more than one time. 
To prove this we will first show that $\varphi$ is $\bigcap$-convex on interval $[t,n].$
Obviously $\varphi (t) >0$. From these facts will follow the
statement of the Theorem~\ref{th33}.

We have
\begin{eqnarray*}
&&\gamma (\ell )=\frac{\sum_{j=0}^{n-\ell}{n-\ell\choose j}\tilde{B}\left( n-\frac{\ell +t}{2}+1-j\right)+
\sum_{j=0}^{n-\ell}{n-\ell\choose j}\tilde{B}\left( n-\frac{\ell +t}{2}-j\right)}{\sum_{j=0}^{n-\ell}{n-\ell\choose j}\tilde{B}\left( n-\frac{\ell +t}{2}-j\right)} \\
&=& 1+\frac{\sum_{j=0}^{n-\ell}{n-\ell\choose j}\tilde{B}\left( n-\frac{\ell +t}{2}+1-j\right)}{\sum_{j=0}^{n-\ell}{n-\ell\choose j}\tilde{B}\left( n-\frac{\ell +t}{2}-j\right)} .
\end{eqnarray*}
Now using identity~(\ref{e111}) we derive the relations
\begin{eqnarray*}
&&\gamma (n,\ell ,t)=\sum_{j=0}^{n-\ell}{n-\ell\choose j}\tilde{B}\left( n-\frac{\ell +t}{2}-j\right) \\
&=& \frac{1}{e}\sum_{ i=1}^\infty \frac{(-1)^{n-\frac{\ell +t}{2}} (i-1)^{n-\frac{\ell +t}{2} }}{i!}\sum_{j=0}^{n-\ell}{n-\ell\choose j}(-1)^j (i-1)^{-j}
=\frac{1}{e}\sum_{i=2}^\infty \frac{(i-1)^{\frac{\ell -t}{2}}(i-2)^{n-\ell}}{i!} .
\end{eqnarray*}
Similar calculations show the validity of the following identity
\begin{eqnarray*}
\gamma (n+2,\ell +2 ,t)&=&\sum_{j=0}^{n-\ell}{n-\ell\choose j}\tilde{B}\left( n-\frac{\ell +t}{2}+1-j\right) 
= \frac{1}{e}\sum_{i=2}^\infty \frac{(i-1)^{\frac{\ell -t}{2}+1}(i-2)^{n-\ell}}{i!}.
\end{eqnarray*}
Hence, for $\gamma (\ell ) -2$,  we have the expression
\begin{eqnarray*}
&&
\gamma (\ell )-2=\frac{\gamma (n+2,\ell +2,t)}{\gamma (n,\ell ,t)} =\frac{ \sum_{i=2}^\infty \frac{(i-1)^{\frac{\ell -t}{2}+1}(i-2)^{n-\ell}}{i!}}{ \sum_{i=2}^\infty \frac{(i-1)^{\frac{\ell -t}{2}}(i-2)^{n-\ell}}{i!}}-1
=\frac{ \sum_{i=2}^\infty \frac{(i-1)^{\frac{\ell -t}{2}}(i-2)^{n-\ell +1}}{i!}}{ \sum_{i=2}^\infty \frac{(i-1)^{\frac{\ell -t}{2}}(i-2)^{n-\ell}}{i!}} .
\end{eqnarray*}
We obtain the following expression for the function $\varphi (\ell ) :$
$$
\varphi (\ell )=t-\ell +2(t-1)\frac{ \sum_{i=2}^\infty \frac{(i-1)^{\frac{\ell -t}{2}}(i-2)^{n-\ell }}{i!}}{  \sum_{i=2}^\infty \frac{(i-1)^{\frac{\ell -t}{2}}(i-2)^{n-\ell +1}}{i!}} .
$$
It is easy to show that the second derivative of this function is negative. This completes the proof of Theorem~\ref{th33} .
  
\bigskip

{\bf II Proof of Theorem~\ref{th21}.}
\bigskip

Denote by $\Omega_0 (n,t) \subset\Omega (n,t)$ the collection of the families of partitions ${\cal A}$ such that $\biggl|\bigcap_{p\in{\cal A}}f(p)\biggr| =0.$
\begin{St}
\label{s1}
\begin{equation}
\label{e122}
\tilde{M}(n,t)=\max_{{\cal A}\in L\tilde{\Omega}(n,t)}|{\cal A}|,
\end{equation}
$$
M_0 (n,t)=\max_{{\cal A}\in\Omega_0 (n,t)}|{\cal A}|=\tilde{M}(n,t).
$$
Moreover, if ${\cal A}\in\tilde{\Omega}(n,t)$ and $\ |{\cal A}|=\tilde{M}(n,t)$, then ${\cal A}\in\Omega_0 (n,t).$
\end{St}
Proof.  First we will prove~(\ref{e122}). For ${\cal A}\in\tilde{\Omega}(n,t)$ assume that $\ |{\cal A}|=\tilde{M}(n,t).$
One can see that either $L(v,w,{\cal A})\in\tilde{\Omega}(n,t)$ or $L(v,w,{\cal A})\in\Omega (n,t)\setminus\tilde{\Omega}(n,t) .$ In the first case we continue shifting. Assume that the second case occurs.
We can assume that $\bigcap_{p\in{\cal A}}f(p) =[t-1]$ and that $ v=t, w=t+1$ and also that  $\cap_{p\in L(v,w,{\cal A})}f(p) =[t].$ Because ${\cal A}$ is maximal, then
\begin{equation}
\label{e2.2}
\left\{ p\in\Omega (n,t):\ [t+1]\subset f(p)\right\} \subset{\cal A}. 
\end{equation}
There are $p_1 ,p_2 \in{\cal A}$ such that
$$
f(p_1)\bigcap [t+1]=[t]
$$
and
$$
f(p_2) \bigcap [t+1]=[t-1]\bigcup\{ t+1\} .
$$
Now we apply the shifting $L(v,w,{\cal A})$ for $v\neq w\in [ n]\setminus\{ t,t+1\}$. We have $\bigcap_{p\in L(v,w,{\cal A})}f(p) =[t-1].$ 

Thus we can assume that $L(v,w,{\cal A})={\cal A}$ for all $v\neq w \in [n]\setminus\{ t,t+1\}$ and
\begin{eqnarray*}
&&f(p_1 )=[a]\setminus\{ t+1\} ,\ a\geq t,\ a\neq t+1,\\
&& f(p_2)=[b]\setminus\{ t\} ,\ b>t .
\end{eqnarray*}
From here and~(\ref{e2.2}) it follows that 
$$
{\cal C}=U(\{[t-1]\bigcup C:\ C\subset [t,\min\{ a,b\}] \})\subset{\cal A} 
$$
and for all $ L(v,w,{\cal C})={\cal C}$ where $v\neq w\in [n]$. Thus $|\bigcap_{p\in{\cal A}}f(p)| <t.$

Now we prove second part of the Statement. Assume that ${\cal A}\subset \tilde{\Omega}(n,t)\setminus\Omega_0 (n,t)$ and $|{\cal A}|=\tilde{M}(n,t).$ We can suppose that ${\cal A}$ is shifted and $\{1\} \in f(p)$ for all $p\in{\cal A} .$
We can also assume that ${\cal A}\in L \tilde{\Omega}(n,t).$ Consider $p\in\Omega (n,t):\ f(p)=\{ 2,\ldots ,n-1 \} .$ Next we will show that $p\in{\cal A},$ which leads to the contradiction of the maximality of ${\cal A}$. Suppose that there exists a partition $p_1 \in {\cal A}$ such that
$$
\biggl| [2,n-1]\bigcap f(p_1 )\biggr|\leq t-1.
$$
We can assume that $f(p_1 )=[t]\bigcup\{n\}
$. We have $p_2 : f(p_2 )=[t-1]\bigcup \{n\}$ belongs to ${\cal A}$ and hence $p_3 : f(p_3 )=[t]$ also belongs to ${\cal A}.$ 
But then $\biggl| f(p_3 )\bigcap f(p_2 )\biggr| =t-1$ which contradicts the $t$-intersecting property of ${\cal A}.$

For further convenience we will make some changes in the definitions, which we will use next.
Let $g({\cal A})$ be the family of subsets of $[n]$ such that ${\cal A}=U(g({\cal A})).$ If ${\cal A}$ is maximal, then we can assume that $g({\cal A})$ is upset and $g^* ({\cal A})$ is the set of its minimal elements. It is easy to see that ${\cal A}\in\Omega (n,t)$ if and only if $g({\cal A})\in I(n,t) $ and $ {\cal A}\in\tilde{\Omega}(n,t)$  if and only if $ g({\cal A})\in\tilde{I}(n,t)$. We can assume that $g({\cal A})$ is left compressed. 

Define
\begin{eqnarray*}&&
s^{+}(a=(a_1 <\ldots <a_j ))=a_j ,\\ &&
s^+ (g({\cal A}))=\max_{a\in g* ({\cal A})} s^{+}(a) ,\\
&&  s_{\min} =\min_{{\cal A}\in L\tilde{\Omega}(n,t):\ |{\cal A}|=\tilde{M}(n,t)}s^{+} (g({\cal A})) .
\end{eqnarray*}
It is easy to see that ${\cal A}\in L\Omega (n,t)$ is a disjoint union
$$
{\cal A} =U_{f\in g^* ({\cal A})} Q(f),
$$
where
$$
Q(f)=\left\{ A\in 2^{[n]} :\ A=f\bigcup B,  B\in [s^{+} (f),n]\right\} ,
$$
and if $f\in g({\cal A})$ is such that $ s^{+}(f)=s^{+}(g({\cal A})) $, then the set of partitions generated only by $f$ is
\begin{equation}
\label{e889}
{\cal A}_f =(U(f)\setminus U(g^* ({\cal A})\setminus\{ f\}))=Q(f) .
\end{equation}
Note also a simple fact that if $f_1 ,f_2 \in g^* ({\cal A})$ and $i\not\in f_1 \bigcup f_2 ,\ j\in f_1 \bigcap f_2 $ for some $
i<j$, then
$$
|f_1 \bigcap f_2 |\geq t+1 .
$$
Next lemma helps us to establish possible sets of $g^* ({\cal A})$  for maximal ${\cal A}\in L\tilde{\Omega} (n,t)$
when $M(n,t)$ is not this maximum. To make the formulation more clear we repeat in Lemma all conditions which we have considered before as default.
\begin{Le}
\label{l4}
For  ${\cal A}\in L\tilde{\Omega} (n,t)$ assume that $ |{\cal A}|=\tilde{M}(n,t)$ and $g({\cal A})\in G({\cal A})$ is such
that $s^{+} (g({\cal A}))=s_{\min} (G({\cal A}))$, then for some $i\geq 2$
$$
g^* ({\cal A})={\cal H}_i .
$$
\end{Le}
Suppose that  $\ell=s^{+}(g({\cal A})) ,\ g_0 ({\cal A})=\{ g\in g^* ({\cal A}):\ s^+ (g)=\ell\}$ and $\ g_1({\cal A}) =g^* ({\cal A})\setminus g_0 ({\cal A}) .$
It is easy to see that $\ell>t+1 .$  From above it follows that if $f_1 ,f_2 \in g_0 ({\cal A})$ and $\biggl| f_1 \bigcap f_2 \biggr| =t$,
then $|f_1 |+|f_2 |=\ell +t .$ 
Denote
$$
\biggl| \bigcap_{f\in g_1 ({\cal A})}f\biggr| =\tau.$$ Consider consequently two cases $\tau <t$ and $\tau \geq t.$

Assume at first that $\tau <t.$ Consider the partition 
$$
g_0 ({\cal A}) =\bigcup_{t<i<\ell} R_i ,\ R_i =g_0 ({\cal A})\bigcap{[n]\choose i} .
$$
Denote $$R^\prime_i =\{ f\subset [\ell -1]:\ f\bigcup \{\ell\}\in R_i \} .
$$
As above, because the set $g({\cal A})$ is left compressed, it follows that for  $$f_i \in R_i^\prime ,\ f_j \in R^\prime_j \  \mbox{and}\ i+j\neq \ell +t,\  |f_i \bigcap f_j |\geq t.$$

Next we show that $R_i =\emptyset .$ 

Assume at first that $\forall R_i \neq\emptyset$ we have $R_{\ell +t-i}=\emptyset $, then for
$$
g^\prime =(g^* ({\cal A})\setminus g_0 ({\cal A}))\bigcup\bigcup_{t<i<\ell}R^\prime_i \in\tilde{I}(n,k)
$$
we have 
$$
|U(g^\prime )|\geq |{\cal A}|\ \hbox{and}\ \ s^{+} (g^\prime )<s^{+} (g({\cal A}))
$$
which contradicts our assumptions.

Now assume that $R_i ,R_{\ell +t-i} \neq \emptyset .$ At first we consider the case when  $i\neq(\ell +t)/2 .$
Consider the new sets
\begin{eqnarray*}
&&\varphi_1 =g_1 ({\cal A})\bigcup \left(g_0 ({\cal A})\setminus \left(R_i \bigcup R_{\ell +t-i}\right)\right)\bigcup R^\prime_i ,\\
&&\varphi_2 =g_1 ({\cal A})\bigcup \left(g_0 ({\cal A})\setminus \left(R_i \bigcup R_{\ell +t-i}\right)\right)\bigcup R^\prime_{\ell +t-i} .
\end{eqnarray*}
We have $\varphi_i \in \tilde{I}(n,k).$ Thus,
$$
{\cal A}_i =U(\varphi_i )\in\tilde{\Omega}(n,t) .
$$
We will show that, under the last assumption,
\begin{equation}
\label{e56}
\max_{j=1,2}|{\cal A}_i |>|{\cal A}|
\end{equation}
and come to a contradiction. Using~(\ref{e889}) it is easy to see that:
\begin{eqnarray*}
&&|{\cal A}\setminus{\cal A}_1| =|R_{\ell +t-i}| \sum_{j=0}^{n-\ell}{n-\ell\choose j}\tilde{B}(n-\ell -t +i-j),\\
&&|{\cal A}_1 \setminus {\cal A}|\geq|R_i |\sum_{j=0}^{n-\ell}{n-\ell\choose j}\tilde{B}(n-i-j+1),\\
&& |{\cal A}\setminus{\cal A}_2 |=|R_i |\sum_{j=0}^{n-\ell}{n-\ell\choose j}\tilde{B}(n-i-j),\\
&&|{\cal A}_2 \setminus {\cal A}| \geq|R_{\ell +t-i}|\sum_{j=0}^{n-\ell}{n-\ell\choose j}\tilde{B}(n-\ell -t+i-j+1) .
\end{eqnarray*}
From these equalities it follows that, if~(\ref{e56}) is not valid, then
$$
|R_{\ell +t-i}| \sum_{j=0}^{n-\ell}{n-\ell\choose j}\tilde{B}(n-\ell -t +i-j)\geq |R_i |\sum_{j=0}^{n-\ell}\tilde{B}(n-i-j+1)
$$
and
$$
|R_i |\sum_{j=0}^{n-\ell}{n-\ell\choose j}\tilde{B}(n-i-j)\geq |R_{\ell +t-i}|\sum_{j=0}^{n-\ell}{n-\ell\choose j}\tilde{B}(n-\ell -t+i-j+1) .
$$
Since $\tilde{B}(n+1)>\tilde{B}(n)$ when $ n>0$, the last two inequalities couldn't be valid together. This contradiction shows that $R_i =\emptyset $ when $ i\neq (\ell +t)/2.$ 

Now consider the case $i=(\ell +t)/2.$ By pigeon- hole principle, there exists $k\in [\ell -1] $ and ${\cal S}\subset R^\prime_{(\ell +t)/2}$ such that $k\not\in B$ for all $\ B\in{\cal S}$ and
\begin{equation}
\label{e9999}
|{\cal S}|\geq\frac{\ell -t}{2(\ell -1)}|R^\prime_{(\ell +t)/2}| .
\end{equation}
Hence, as before, we have $\biggl| B_1 \bigcap B_2 \biggr|\geq t$ for all $\ B_1 ,B_2 \in{\cal S}$ and
$$
f^\prime =(g^* ({\cal A})\setminus R_{(\ell +t)/2})\bigcup {\cal S}\in\tilde{I}(n,t).
$$
Next we show that 
\begin{equation}
\label{e4.9}
|U(f^\prime )|> |{\cal A}|.
\end{equation}
Consider the partition
$$
{\cal A}={\cal G}_1 \bigcup{\cal G}_2 ,
$$
where
\begin{eqnarray*}
&&
{\cal G}_1 =U(g^* ({\cal A})\setminus R_{(\ell +t)/2}),\\
&&
{\cal G}_2 =U(R_{(\ell +t)/2})\setminus U(g^* ({\cal A})\setminus R_{(\ell +t)/2}) .
\end{eqnarray*}
Consider also the partition
$$
U(f^\prime )={\cal G}_1 \bigcup{\cal G}_3 ,
$$
where
$$
{\cal G}_3 =U({\cal S})\setminus U(g^* ({\cal A})\setminus R_{(\ell +t)/2}) .
$$
We should show that
\begin{equation}
\label{e4.10}
|{\cal G}_3 |>|{\cal G}_2 |.
\end{equation}
We have
$$
|{\cal G}_2 |=|R_{(\ell +t)/2}|\sum_{j=0}^{n-\ell}{n-\ell\choose j}\tilde{B}\left( n-\frac{\ell +t}{2}-j\right)
$$
and
$$
|{\cal G}_3 |\geq |{\cal S}|\sum_{j=0}^{n-\ell +1}{n-\ell +1 \choose j}\tilde{B}\left( n-\frac{\ell +t}{2}-j+1\right) .
$$
Hence, for~(\ref{e4.10}) to be true, it is sufficient that
$$
\frac{\ell -t}{2(\ell -1)}\gamma (\ell )>1.
$$
The last inequality is true because, otherwise,
from~(\ref{th21}) it follows that  $\tilde{M}(n,k)=M(n,k).$ Hence $R_{\frac{\ell +t}{2}}=\emptyset .$

Now consider the case $\tau\geq t.$ We have
$$
\bigcap_{f\in g_1 ({\cal A})} f=[\tau ],
$$
$$
\ell =s^{+} (g({\cal A}))>\tau 
$$
and for all $f\in g_0 ({\cal A})$,
\begin{eqnarray*}
&& \biggl| F\bigcap [\tau ]\biggr|\geq\tau -1,\\
&& \hbox{if}\  |f\cap [\tau ]|=\tau -1 , \hbox{then}\  [\tau +1,\ell ]\in f.
\end{eqnarray*}
Let's show that $\tau \leq t+1.$

If $\tau\geq t+2$, then, for $f_1 ,f_2 \in g({\cal A}),$
$$
\biggl| f_1 \bigcap f_2 \bigcap[\tau ]\biggr|\geq\tau -2\geq t
$$
and thus, setting $g^\prime_0 ({\cal A})=\{ f\subset [\ell -1]: f\bigcup\{\ell\}\in g_0 ({\cal A})\} ,$
we have
$$
\varphi =(g^* ({\cal A})\setminus g_0 ({\cal A}))\bigcup g^\prime_0 ({\cal A})\in\tilde{I}(n,k)
$$
and
\begin{eqnarray*}
&&|U(\varphi )|\geq|{\cal A}|,\
 s^{+}(\varphi )<\ell .
\end{eqnarray*}
This gives us the contradiction of minimality of $\ell .$

Assume now that $\tau =t+1.$ In this case it is necessary that $\ell =t+2$. Otherwise, using the argument above (deleting $\ell $ from each element of $g_0 ({\cal A})$), we end up generating the set $\varphi \in\tilde{I}(n,k)$ for which $|U(\varphi )|\geq |{\cal A}|$ and $s^{+}(\varphi )<\ell .$ It is clear that $\tau =t+1$ and $\ell =t+2$, then $g^* ({\cal A})={\cal H}_2 .$

At last, consider the case $\tau =t.$ Define $g^\prime_0 ({\cal A})=\left\{ f\in g_0 ({\cal A}): | f\bigcap [t]| =t-1\right\} .$ 
We have
$$
g^\prime_0 ({\cal A})\subset\{ f\subset [\ell ]:\ |f\cap [t]|=t-1, [t+1,\ell ]\subset f\} 
$$
and for $f\in g^* ({\cal A})\setminus g_0^\prime ({\cal A})$ we have $[t]\subset f$ and $| f\bigcap [t+1,\ell ]|\geq 1.$ 

Hence 
$$
U(g^* ({\cal A})) \subset U({\cal H}_{\ell -t} ).
$$
Since ${\cal A}$ is maximal, $g^* ({\cal A})={\cal H}_{\ell -t} .$ 
Family ${\cal H}_{n-t}$ is trivially $t$-intersecting, so we can assume that $i<n-t.$
Denote $S_i =|U({\cal H}_i )|$. Next we will prove that if 
$S_i <S_{i+1}$, then $S_{i+1}<S_{i+2} .$
We have
$$
S_i =(n-i)!-\sum_{j=0}^{n-t-i}{n-t-i\choose j}\tilde{B}(n-t-j)+t\sum_{j=0}^{n-t-i}\tilde{B}(n-t-i-j+1) 
$$
and we should show that from inequality
\begin{equation}
\label{e990}
\sum_{j=0}^{n-t-i-1}{n-t-i-1\choose j}\tilde{B}(n-t-j+1)\geq t\sum_{j=0}^{n-t-i-1}{n-t-i-1\choose j}\tilde{B}(n-t-j-i+1)
\end{equation}
follows
\begin{eqnarray}
\label{e91}
&&\sum_{j=0}^{n-t-i-2}{n-t-i-2\choose j}\tilde{B}(n-t-j+1)
\geq t\sum_{j=0}^{n-t-i-2}{n-t-i-2\choose j}\tilde{B}(n-t-j-i).
\end{eqnarray}
We rewrite inequality~(\ref{e990}) as follows
\begin{eqnarray*}
&&\sum_{j=0}^{n-t-i-2}{n-t-i-2\choose j}\tilde{B}(n-t-j+1)
+\sum_{j=0}^{n-t-i-2}{n-t-i-2\choose j}\tilde{B}(n-t-j)\\
&\geq& t\sum_{j=0}^{n-t-i-2}{n-t-i-2\choose j}\tilde{B}(n-t-i-j+1)
+t\sum_{j=0}^{n-t-i-2}{n-t-i-2\choose j}\tilde{B}(n-t-i-j) .
\end{eqnarray*}
From here, it is clear that if~(\ref{e91}) is true, then~(\ref{e90}) is also true.
From here and expressions for $S_2$ and $S_{n-t-1}$ follows the statement of Theorem~\ref{th21}.
Since, for fixed $t$,
$$
\frac{\sum_{j=1}^{n-t-2}{n-t-2\choose j}\tilde{B}\left(n-t-j\right)}{\sum_{j=0}^{n-t-2}{n-t-2\choose j}\tilde{B}\left( n-t-1-j\right)}\to\infty ,\ n\to\infty ,
$$
and
$$
\frac{\tilde{B}(n-t-1)}{\sum_{j=1}^{n-t-2}{n-t-2\choose j}\tilde{B}\left(n-t-j\right)}\to 0,\ n\to\infty ,
$$
it 
 follows that
for sufficiently large $n$ and fixed $t$:
\begin{eqnarray*}
&&
S_2= B(n-t)-\tilde{B}(n-t)-\tilde{B}(n-t-1)+t >S_{n-t-1}=B(n-t)\\
&-&\sum_{j=0}^{n-t-2}{n-t-2\choose j}\tilde{B}(n-t-j)+t\sum_{j=0}^{n-t-2}{n-t-2\choose j}\tilde{B}(n-t-j-1).
\end{eqnarray*}
Therefore, for $n>n_2 (t),$ 
$$
\tilde{ M}(n,t)=B(n-t)-\tilde{B}(n-1)-\tilde{B}(n-t-1)+t .
$$
\newpage
\begin{center}
{\bf Appendix}
\end{center}\bigskip
Let us remark that FKG inequality says that for $\mu : 2^{[m]} \to R_+$  such that
 \begin{equation}
 \label{e44}
 \mu (a)\mu (b) \leq \mu \left( a\bigcap b\right)\mu \left(a\bigcup b\right) ,\ a,b\in 2^{[m]} ,
 \end{equation}
 and  for a pair of nondecreasing functions $f_1 ,f_2 : 2^{[m]}\to R$, the following 
 inequality is valid:
 \begin{equation}
 \label{errr}
 \sum_{Y\in 2^{[m]}}\mu (Y)f_1 (Y)\sum_{Y\in 2^{[m]}}\mu (Y)f_2 (Y)\leq\sum_{Y\in 2^{[m]}}\mu (Y)f_1 (Y)f_2 (Y)\sum_{Y\in 2^{[m]}}\mu (Y) .
 \end{equation}
 Now we choose
 $$
 \mu (Y)=\tilde{B}\left( n-\frac{\ell +t}{2} -|Y|\right) .
 $$
 Note that if~(\ref{e44}) is true for this choice of $\mu$, then setting $f_1 =I_{X\in\Gamma :\ x\in X}$ and
 $f_2 =I_{X\in 2^{[n-(\ell +t)/2-1]},\ x\in X}$
 in~(\ref{errr}) proves inequality~(\ref{e788}).   

 Define $\bar{a}=n-t-r-a,\ \bar{b}=n-t-r-b,\ \bar{\delta}=n-r-t-\delta$. Then inequality~(\ref{e44}) is equivalent to inequality
 $$
 \tilde{B}(\bar{a})\tilde{B}(\bar{b})\leq\tilde{B}(\bar{a}\bigcap\bar{b})\tilde{B}(\bar{a}\bigcup\bar{b}).
 $$
 Function $F(i) \geq 0,\ i\in Z_+$ is called log - convex if
 \begin{equation}
 \label{e22}
 F(i+1)F(i-1)\geq F^2(i).
 \end{equation}
 FKG - condition~(\ref{e44}) is equivalent to the log - convexity property of $\tilde{B}.$ We are going to demonstrate   that $F(|Y|)=\mu (Y)=\tilde{B}\left( n-\frac{\ell +t}{2} -|Y|\right)$ 
 satisfy inequality~(\ref{e22}) for all possible $|Y|,$ except $n-t -r-|Y|\neq 2,4.$ 
  \begin{Le}\label{p1}
 Inequality
 \begin{equation}
 \label{e23}
 \tilde{B}(k+1)\tilde{B}(k-1)\geq \tilde{B}^2(k)
 \end{equation}
 is true for $k\in Z_+ \setminus \{2,4\}$
 \end{Le}
 From above considerations it  follows it is left to consider the case $k=2m,\ m>2$.  We will prove this lemma by using asymptotic of $\tilde{B}(k)$. 
 
 Next part of text we devote to finding the asymptotic fo $\tilde{B}(n)$.
 
 \begin{Le} \label{p2}
 The following asymptotic of $\tilde{B}(n)$ is true
 \begin{equation}
 \label{d1}
 \tilde{B}(n)=\frac{n! \exp(e^r-r-1)}{r^{n}(4\pi B)^{1/2}}(1\pm 11e^{-r} ),\ r\geq 12.
 \end{equation}
 where r satisfy equality $r(e^r-1)=n$ and $B=\frac{1}{2}r((r+1)e^r -1)$.
  \end{Le}
 {\bf Proof}
 
 To proof this lemma we will use The Moser - Wyman expansion of the Bell numbers~\cite{12}. We will follow the text~\cite{13} and introduce the extension of proof from~\cite{13} for extend Bell number
 $\tilde{B}(n)$ for completeness (it is quite similar, besides we need calculate the explicit bounds for rest term of asymptotic also).
 
 Because
 $$
 \sum_{n=0}^\infty \frac{\tilde{B}(n)x^n}{n!}=\exp(e^x -x-1),  
 $$
 using  Cauchy's formula, we have
 $$
 \frac{\tilde{B}(n)}{n!} =\frac{1}{2\pi i}\oint_{|z|=r} \frac{\exp(e^z -z-1)}{z^{n+1}}dz.
 $$
 Contour integration yields 
 $$
 \tilde{B}(n)=\frac{n!}{2\pi r^n}\int_{-\pi}^\pi \exp (e^{re^{i\theta}} -re^{i\theta} -in\theta -1)d\theta .
 $$
 Define
 $$
 F(\theta )=e^{re^{i\theta}} - re^{i\theta} -in\theta -(e^r -r).
 $$
 We have
 $$
 \tilde{B}(n)=A\int_{-\pi}^\pi \exp(F(\theta ))d\theta 
 $$
 where
 $
 A=\frac{n! \exp (e^r -r-1)}{2\pi r^n} .
 $
 
 Define $$
 \epsilon =e^{-\frac{3}{8}r},\ J_1=\int_{-\pi}^{\epsilon}\exp(F(\theta ))d\theta,\ J_2 =\int_{\epsilon}^{\pi}\exp(F(\theta))d\theta .
 $$
 Using inequality
 $\cos(\theta) \leq 1-\frac{\epsilon^2}{2}+\frac{\epsilon^4}{24}$, we have
 \begin{eqnarray*}&&
 \tilde{B}(n) =AJ_1 +AJ_2 +A\int_{-\epsilon}^{\epsilon} \exp (F(\theta))d\theta , \\
 && |\exp(F(\theta))|=\exp(Re (F(\theta))) =\exp\left( e^{r\cos(\theta)}\cos(r\sin(\theta)) -e^r +r(1-\cos (\theta ))\right)\\
 && \leq \exp\left( e^{r\cos(\theta)} -e^r +r(1-\cos (\theta ))\right) = \exp\left( e^r\left(e^{r(\cos(\theta)-1)} - 1\right) +r(1-\cos (\theta ))\right)\\
 && \leq\exp\left( e^r \left(e^{r\left(1-\frac{\theta^2}{2}+\frac{\theta^4}{24}-1\right)}-1\right)+r\right) = \exp\left( e^r \left(e^{-\frac{\theta^2}{2}\left(1-\frac{\theta^2}{12}\right)r} -1\right)+r\right)\\
 && \leq\exp\left( e^r \left(  -\frac{\epsilon^2}{2}\left(1-\frac{\epsilon^2}{12}\right)r+\frac{\epsilon^4 \left(1-\frac{\epsilon^2}{12}\right)^2}{8}r^2\right)+r\right)\\
&&  \leq \exp\left(-\frac{1}{2}e^{r}r\epsilon^2 \left( 1-\frac{\epsilon^2}{12}\right) \left(1-\frac{\epsilon^2\left(1-\frac{\epsilon^2}{12}\right)r}{4}\right) +r\right) 
 \\
 && \leq \exp\left(-\frac{1}{2}e^{r/4}r \left( 1-\frac{e^{-3r/4}}{12}\right) \left(1-\frac{e^{-3r/4}\left(1-\frac{e^{-3r/4}}{12}\right)r}{4}\right) +r\right) .
 \end{eqnarray*}
Because
$$
 1-\frac{e^{-3r/4}}{12} >\frac{11}{12},\  1-\frac{e^{-3r/4}\left(1-\frac{e^{-3r/4}}{12}\right)r}{4}>\frac{3}{4},\ \hbox{when}\ r>12,
 $$
 we have
 $$
 |\exp(F(\theta))| \leq\exp\left(-r\frac{1}{3} e^{r/4}+r\right)< e^{-\frac{r}{4}e^{r/4} }
 $$
 and, hence
 $$
 J_2 \leq A\pi e^{-\frac{r}{4}e^{r/4}} .
 $$
 We have
 $$\tilde{B}(n)=A\left(\int_{-\epsilon}^{\epsilon}\exp(F(\theta))d\theta \pm e^{-\frac{r}{4}e^{r/4}} \pi\right).
 $$
 Consider the expansion
 $$
 F(\theta )=(re^r -r-n)i\theta -\frac{1}{2} (r^2 e^r +re^r -r)\theta^2 +\sum_{k=3}^\infty \left(\left (r\frac{d}{r}\right)^k (e^r) -r\right) (i\theta)^k,\ r(e^r-1)=n.
 $$
 Hence we have
 $$
 F(\theta )=-\frac{1}{2}(r^2 e^r +r e^r -r)\theta^2 + \sum_{k=3}^\infty \left(\left (r\frac{d}{r}\right)^k (e^r) -r\right) (i\theta)^k.
 $$
 Define
 $$
 \phi =\left(\frac{1}{2}(r^2e^r +re^r -r)\right)^{1/2} \theta,\ a_k =\frac{\left(e^{-r}\left(r\frac{d}{dr}\right)^{k+2}(e^r) -re^{-r}\right)(i\phi)^{k+2}}{(k+2)!\left(\frac{1}{2}(r^2 +r-re^{-r}\right)^{\frac{k+2}{2}}}, z=e^{-r/2},
 $$
 $$
 f(z)=\sum_{k=1}^\infty a_k z^k.
 $$
 Then 
 \begin{eqnarray*}
 &&F(\theta )=-\phi^2 +f(z),\ \bar{B}(n)=C\left(\int_{-h}^h \exp (-\phi^2 +F(z))dz \pm \pi e^{-\frac{r}{4}e^{r/4}}\left(\frac{1}{2} r(r+1)e^{r}-r\right)^{1/2}\right),\\
 && h=\left(\frac{1}{2}r((r+1)e^r -1)\right)^{1/2}e^{-3r/8},\ C=\frac{A}{ \left(\frac{1}{2}r((r+1)e^r -1)\right)^{1/2}} .
 \end{eqnarray*}
 Consider the expansion
 $$
 e^{f(z)}=\sum_{k=0}^\infty b_k z^k,\ b_0 =e^{f(0)} =1,\ b_1 =e^{f(0)}f^\prime (0)=a_1,\ b_2 =a_2+\frac{a_1^2}{2}.
 $$
 We have
 \begin{eqnarray*}
 &&
 |a_k|=\Biggl| \frac{\left(\sum_{m=1}^{k+2}S(k+2,m)r^m -re^{-r}\right)(i\phi)^{k+2}}{(k+2)!\left(\frac{1}{2}(r^2+r -re^{-r})\right)^{(k+2)/2}}\Biggr|\\
 &&\leq \frac{2^{(k+2)/2} |\phi|^{k+2} B(k+2)}{(k+2)!} < |2\phi|^{k+2} .
 \end{eqnarray*}
 Here $S(m,k)$ is Stirling number of second kind, $B(n)$ is Bell number. We used inequalities $B(n)\leq n!,\ \left(r\frac{d}{dr}\right)^{k+2}(e^r) =\sum_{m=1}^{k+2}S(k+2,m)r^n e^r\leq r^{k+2}B(k+2)e^r.$ 
 We use formula for coefficients $b_k$ in composite function $\exp \left(\sum_{k=3}^\infty a_k z^k\right)=\sum_{k=1}^\infty b_k z^k):$
 \begin{eqnarray*} &&
 b_m =\sum_{\sum_{p=1}^n pj_p =m,\ \sum_{p=1}^m j_p =k} \frac{1}{\prod_{p=1}^n j_p!}\prod_{p=1}^m a_p^{j_p}\\
 && \leq \sum_{\sum_{p=1}^m pj_p =m,\ \sum_{p=1}^m j_p =k} \frac{1}{\prod_{p=1}^m j_p!}\prod_{p=1}^m (2\phi)^{(p+2)j_p}\leq (2\phi)^{m} \sum_{ \sum_{p=1}^m pj_p =m,\ \sum_{p=1}^m j_p =k} \frac{(2\phi)^{2k}}{\prod_{p=1}^m j_p!} \\&&
 =(2\phi)^{m}\sum_{\sum_{p=1}^m pj_p =m}  \sum_{\sum_{p=1}^m j_p =k} \frac{(2\phi)^{2k}}{\prod_{p=1}^m( j_p)!} \
 =(2\phi)^{m}\sum_{\sum_{p=1}^m pj_p =m} \frac{(2\phi)^{2k}}{k!} \sum_{\sum_{p=1}^m j_p =k} \frac{k!}{\prod_{p=1}^m (j_p)!} \\ 
 &&=(2\phi)^{m}\sum_{\sum_{p=1}^m pj_p =m} \frac{(2\phi)^{2k}}{k!} {m-1\choose k-1}\leq  (2\phi)^{m}\sum_{\sum_{p=1}^m pj_p =m} \frac{(2\phi)^{2k}}{k!} {m-1\choose k-1}\\ &&
  \leq (2\phi )^m (1+(2\phi)^2 )^{m-1} .
 \end{eqnarray*}
 Next we have
 \begin{eqnarray*}&&
 \Biggl|\sum_{k=s}^\infty b_k z^k\Biggr| \leq \left(|(2\phi)|^{s+2}(1+|2\phi|^2 )^{s-1}|z|^s \right)\sum_{i=0}^\infty \mu^i =\frac{|2\phi|^{s+2}(1+|2\phi|^2)^{s-1}|z|^s}{M},\\ && M=1-|z||2\phi|(1+|2\phi|^2)\ \mu =|2\phi|(1+|2\phi|^2)|z|<1 .
 \end{eqnarray*}
 We impose conditions $M>\frac{1}{2},\ |\phi|\leq h , z=e^{-r/2}.$
 We have
 \begin{eqnarray}&&
 \label{j1}
 \tilde{B}(n) =C\biggl(\int_{-h}^h e^{-\phi^2} d\phi + e^{-r}\int_{-h}^h b_2 e^{-\phi^2} d\phi \pm\int_{-\infty}^\infty \Biggl|\sum_{k=2}^\infty b_{2k}(\phi)e^{-kr}\Biggr|e^{-\phi^2} d\theta\\ \nonumber
 && \pm \pi e^{-\frac{r}{4}e^{r/4}}\left(\frac{1}{2} (r(r+1)e^{r}-r)\right)^{1/2}\biggr) 
 \end{eqnarray}
 Because 
 \begin{equation}
 \label{r5}
 \int_{h}^\infty e^{-\phi^2}d\phi =\frac{1}{2h}e^{-h^2} \left(1\pm\frac{1}{h^2}\right),\end{equation}we have
 $$
 \int_{-h}^h e^{-\phi^2} d\phi  = \int_{-\infty}^{\infty} e^{-\phi^2} d\phi - 2\int_{h}^\infty e^{-\phi^2} d\phi = \sqrt{\pi} -\frac{1}{h}e^{-h^2}\left(1\pm \frac{1}{h^2}\right) .
 $$
 It follows the asymptotic equation
 \begin{eqnarray}&&\label{f5}
 \bar{B}(n) =C\biggl(e^{-r}\int_{-h}^h b_2 e^{-\phi^2} d\phi  + \sqrt{\pi} -\frac{1}{h}e^{-h^2}\left(1\pm \frac{1}{h^2}\right) \pm e^{-2r} \int_{-\infty}^\infty  (2\phi)^{4}(1+(2\phi)^2)^{3}e^{-\phi^2} d\theta\\ \nonumber
 && \pm \pi e^{-\frac{r}{4}e^{r/4}}\left(\frac{1}{2} (r(r+1)e^{r}-r)\right)^{1/2}\biggr) \nonumber
 \end{eqnarray}
 Next we have
 \begin{eqnarray*}
 && a_1 =\sqrt{2} \frac{r^3+3r^2+r-re^{-r}}{3(r^2 +r-re^{-r})^{3/2}} (i\phi)^3;\\
 && a_2 =\frac{r^4 +6r^3 +7r^2 +r-re^{-r}}{6(r^2 +r-re^{-r})^2} \phi^4;\\
 &&b_2 =a_2+\frac{a_1^2}{2}= \frac{r^4 +6r^3 +7r^2 +r-re^{-r}}{6(r^2 +r-re^{-r})^2} \phi^4 - \frac{(r^3+3r^2+r-re^{-r})^2 }{9(r^2 +r-re^{-r})^3} \phi^6.
 \end{eqnarray*}
Integrating in parts and using asymptotic~(\ref{r5}) we have
\begin{eqnarray*}
&&\int_{-h}^h \phi^4 e^{-\phi^2}d\phi=-h^3 e^{-h^2} -\frac{3}{2}he^{-h^2} +\frac{3}{4}\left(\sqrt{\pi}-2\int_h^\infty e^{-\phi^2}d\phi\right)=-h^3 e^{-h^2} -\frac{3}{2}he^{-h^2} \\&& +\frac{3}{4}\left(\sqrt{\pi}-2\int_h^\infty e^{-\phi^2}d\phi\right)
=  -h^3 e^{-h^2} -\frac{3}{2}he^{-h^2} +\frac{3}{4}\left(\sqrt{\pi}-\frac{1}{h}e^{-h^2}\left(1\pm \frac{1}{h^2}\right)\right),\\ 
&&\int_{-h}^h \phi^6 e^{-\phi^2}d\phi =-h^5e^{-h^2}+\frac{5}{2}\int_{-h}^h \phi^4 e^{-\phi^2}d\phi\int_{-h}^h \phi^8 e^{-\phi^2}d\phi =-h^7 e^{-h^2}+\frac{7}{2}\int_{-h}^h \phi^6 e^{-\phi^2}d\phi  \int_{-h}^h \phi^{10} e^{-\phi^2}d\phi \\
&&=-h^9e^{-h^2}+\frac{9}{2}\int_{-h}^h \phi^{8} e^{-\phi^2}d\phi.
\end{eqnarray*}
Because $\frac{re^{r/8}}{\sqrt{2}}< h< re^{r/8}\sqrt{2},\ r^8 e^r > 100 r^6 e^{3r/4}, $ we have

\begin{eqnarray*}&&
\int_{-h}^h (2\phi)^4 (1+(2\phi)^2)^3 e^{-\phi^2}d\phi =16\int_{-\infty}^\infty  (\theta^4 +12\theta^6 +48\theta^8 +64\theta^{10} )e^{-\phi^2}d\phi\\
&&= -32\left(e^{-h^2}h\left(\frac{5133}{2}+1711h^2+334h^4+129h^6+32h^8\right) -\frac{5133}{3}\left(\sqrt{\pi}-\frac{e^{-h^2}}{h}\left( 1\pm\frac{1}{h^2}\right)\right)\right)\\
&& 
=-32\left(e^{-\frac{r^2}{2} e^{r/4}}( 640r^8 e^r +1200r^6 e^{3r/4}+1600r^4 e^{r/2}+4000re^{r/4}+3000)+5133\left(\sqrt{\pi} + \frac{e^{-\frac{r^2}{2}e^{r/4}}}{e^{r/2}r}\right)\right)\\
&&= \pm 2^{17} \left( e^{-\frac{r^2}{2}e^{r/4}} r^8 e^r+1\right),
\end{eqnarray*}
\begin{eqnarray*}
&& \int_{-h}^h b_2 e^{-\phi^2}d\phi =\int_{-h}^h a_2 e^{-\phi^2}d\phi +\frac{1}{2}\int_{-h}^h a_1^2 e^{-\phi^2}d\phi \\
&& =\frac{r^4 +6r^3 +7r^2 +r-re^{-r}}{6(r^2 +r-re^{-r})^2} \left(  -h^3 e^{-h^2} -\frac{3}{2}he^{-h^2} +\frac{3}{4}\left(\sqrt{\pi}-\frac{1}{h}e^{-h^2}\left(1\pm \frac{1}{h^2}\right)\right)\right)\\
&&-\frac{1}{2} \frac{(r^3+3r^2+r-re^{-r})^2 }{9(r^2 +r-re^{-r})^3} \left(-h^5 e^{-h^2} -\frac{5}{2}h^3 e^{-h^2} -\frac{15}{4}he^{-h^2} +\frac{15}{8}\left(\sqrt{\pi}-\frac{1}{h}e^{-h^2}\left(1\pm\frac{1}{h^2}\right)\right)\right)\\
&&=\left(\frac{\sqrt{\pi}}{8} \frac{r^4 +6r^3 +7r^2 +r-re^{-r}}{6(r^2 +r-re^{-r})^2} -\frac{5\sqrt{\pi}}{24} \frac{(r^3+3r^2+r-re^{-r})^2 }{9(r^2 +r-re^{-r})^3}\right)\\
&&-e^{-h^2} \frac{r^4 +6r^3 +7r^2 +r-re^{-r}}{6(r^2 +r-re^{-r})^2}\left(h^3+\frac{3}{2}h+\frac{3}{4h}\right)\\
&& +\frac{1}{2} e^{-h^2}\frac{(r^3+3r^2+r-re^{-r})^2 }{9(r^2 +r-re^{-r})^3}\left(h^5 +\frac{5}{2}h^3 +\frac{15}{4}h \frac{15}{8h}\right) \\
&&\pm e^{-h^2}\left(\frac{r^4 +6r^3 +7r^2 +r-re^{-r}}{6(r^2 +r-re^{-r})^2} + \frac{1}{2} \frac{(r^3+3r^2+r-re^{-r})^2 }{9(r^2 +r-re^{-r})^3}\right) \frac{1}{h^3}\\ \nonumber
&& = \left(\frac{\sqrt{\pi}}{8} \frac{r^4 +6r^3 +7r^2 +r-re^{-r}}{6(r^2 +r-re^{-r})^2} -\frac{5\sqrt{\pi}}{24} \frac{(r^3+3r^2+r-re^{-r})^2 }{9(r^2 +r-re^{-r})^3}\right)\\
&&\pm e^{-\frac{r^2}{2}e^{r/4}}\left( \frac{r^4 +6r^3 +7r^2 +r-re^{-r}}{6(r^2 +r-re^{-r})^2}4r^3 e^{3r/8} +\frac{1}{2}\frac{(r^3+3r^2+r-re^{-r})^2 }{9(r^2 +r-re^{-r})^3} r^5 e^{5r/8}\right)
\end{eqnarray*}
We are ready to write the asymptotic
\begin{eqnarray}
\nonumber&&
\bar{B}(n) =C\sqrt{\pi}\biggl(1+\frac{1}{\sqrt{\pi}}e^{-r}\int_{-h}^h b_2 e^{-\phi^2} d\phi   \pm \frac{1}{\sqrt{\pi}}\biggl(2\sqrt{\pi}e^{-\frac{r^2}{2}e^{r/4}}+e^{-2r} \int_{-\infty}^\infty  (2\phi)^{4}(1+(2\phi)^2)^{3}e^{-\phi^2} d\theta\\&& + \pi e^{-\frac{r}{4}e^{r/4}}\left(\frac{1}{2} (r(r+1)e^{r}-r)\right)^{1/2}\biggr)\biggr)\\ \label{p2}
&&\nonumber=C\sqrt{\pi}\left(1+\frac{1}{\sqrt{\pi}}e^{-r}\int_{-h}^h b_2 e^{-\frac{r^2}{2}e^{r/4}} d\phi   \pm \frac{1}{\sqrt{\pi}}\left(2\sqrt{\pi}e^{-\frac{r^2}{2}e^{r/2}} +2^{17} \left( e^{-\frac{r^2 e^{r/4}}2} r^8 e^r+1\right)e^{-2r} +\sqrt{\pi} e^{-\frac{r}{4}e^{r/4}}re^{r/2}\right)\right)\\ \nonumber
&& =C\sqrt{\pi}\Biggl(1+ e^{-r} \left(\frac{1}{8} \frac{r^4 +6r^3 +7r^2 +r-re^{-r}}{6(r^2 +r-re^{-r})^2} -\frac{5}{24} \frac{(r^3+3r^2+r-re^{-r})^2 }{9(r^2 +r-re^{-r})^3}\right)\\ \nonumber
&&\pm \frac{1}{\sqrt{\pi}}\Biggl(e^{-\frac{r^2}{2} e^{r/4}}\biggl( \frac{r^4 +6r^3 +7r^2 +r-re^{-r}}{6(r^2 +r-re^{-r})^2}4r^3 e^{-5r/8} +\frac{1}{2}\frac{(r^3+3r^2+r-re^{-r})^2 }{9(r^2 +r-re^{-r})^3} r^5 e^{-3r/8}\\&&
+2^{18}r^8e^{-r}+2\pi\biggr)+2re^{-\frac{r}{4}e^{r/4} +\frac{r}{2}}\Biggr)\Biggr)= C\sqrt{\pi}( 1\pm 11e^{-r}),\ r>10. \nonumber
 \end{eqnarray}
   Through formula~(\ref{e111}) it is easy to see that  when $a,b,\delta,a+b-\delta <n-t-r$, the inequality~(\ref{e33}) follows from the inequality
 \begin{eqnarray}&&\label{f1}
 \left((-1)^{n-\frac{\ell +t}{2}-a}+\frac{1}{2}+\sum_{k=2}^\infty \frac{k^{n-\frac{\ell+t}{2}-a}}{(k+1)!}\right)\left((-1)^{n-\frac{\ell +t}{2}-b}+\frac{1}{2}+\sum_{k=2}^\infty \frac{k^{n-\frac{\ell+t}{2}-b}}{(k+1)!}\right)\\ \nonumber
 &&\leq\left((-1)^{n-\frac{\ell +t}{2}-\delta}+\frac{1}{2}+\sum_{k=2}^\infty \frac{k^{n-\frac{\ell+t}{2}-\delta}}{(k+1)!}\right)\left((-1)^{n-\frac{\ell +t}{2}-(a+b)+\delta}+\frac{1}{2}+\sum_{k=2}^\infty \frac{k^{n-\frac{\ell+t}{2}-(a+b)+\delta}}{(k+1)!}\right)
 \end{eqnarray}
 Define $ \bar{a}= n-r-t-a,\ \bar{b}=n-r-t-b,\ \bar{\delta} =n-r-t-\delta =\bar{a}\bigcap\bar{b},\ n-r-t-(a+b)+\delta =\bar{a}\bigcup\bar{b}$ Inequality~(\ref{f1}) is equivalent to the inequality
 \begin{equation}
 \label{d1}
 \tilde{B}(\bar{a})\tilde{B}(\bar{b}) \leq \tilde{B}(\bar{a}\bigcap\bar{b})\tilde{B}(\bar{a}\bigcup\bar{b}).
 \end{equation}  
 As we noticed before last inequality equivalent to log - convexity of function $\bar{B}(n)$:
 \begin{equation}\label{l2}
 \tilde{B}^2 (k)\leq \tilde{B}(k+1)\tilde{B}(k-1) ,\ k\in N.
 \end{equation}
Next, using asymptotic~(\ref{d1}) we prove inequality~(\ref{e23}) for sufficiently large $n>n_0$ where $n_0<e^{12}$. Then using software "Mathematica" we show that inequality~(\ref{l2}) is true for all other values of $n=2m <e^{12},\ m>2$. This complete the proof of Lemma~\ref{p1}.
 Simple calculations show the validness of the inequality
 $$
 D(n)\stackrel{\Delta}{=}\sqrt{\pi}C(n)=\sqrt{\pi}C=\frac{n!e^{\frac{n}{r}}}{\sqrt{2\pi} r^{n-1}(n+r) (r^2+n(r+1))^{1/2}} .
 $$
 Following inequality is valid:
 $$
 r(n+1)<r(n)+\frac{1}{n}.
 $$
 Indeed
 $$
 \left(r(n)+\frac{1}{n}\right)\left(e^{r+\frac{1}{n}} -1\right) >n+1.
 $$
 Hence for $r=r(n),\ \ln (n) >r(n)>10$, we have
 \begin{eqnarray*}
 &&
 C(n-1)C(n+1) -(C(n))^2 >\frac{(n-1)! (n+1)! e^{\frac{n-1}{r} +\frac{n+1}{r+\frac{1}{n}}}}{2\pi r^{n-2} \left(r+\frac{1}{n}\right)^{n}\left(n+r+\frac{1}{n}\right)}\\
 &&\cdot \frac{1}{(n-1+r)\left(r^2 +(n-1)(r +1)\right)^{1/2}\left(\left(r+\frac{1}{n}\right)^2 +(n+1)\left(r+1+\frac{1}{n}\right)\right)^{1/2}} \\
 && -\frac{n!^2 e^{2\frac{n}{r}}}{2\pi r^{2(n-1)}(r+n)^2 (r^2+n(r+1))}> \left( \frac{e^{r^{-2}}}{(1-1/n)} -1\right) \frac{n!^2 e^{2\frac{n}{r}}}{2\pi r^{2n} (r+n)^2 (r^2+n(r+1))}.
 \end{eqnarray*}  
 Using inequalities $\frac{n}{n-1}>1+\frac{1}{n},\ e^{r^{-2}} > 1+r^{-2}$ and, hence $ \frac{e^{r^{-2}}}{(1-1/n)}>1+r^{-2}$ we have
 $$
  C(n-1)C(n+1) -(C(n))^2 >  V\stackrel{\Delta}{=}\frac{n!^2 e^{2\frac{n}{r}}}{2\pi r^{2n+2} (r+n)^2 (r^2+n(r+1))}.
  $$
  Let's $\tilde{B}(n) =\sqrt{\pi}C(n)(1+d(n))$. We need to prove the inequality
  $$
  C(n-1)(1+d(n-1))C(n+1)(1+d(n+1)) \geq C^2 (n)(1+d(n))^2
  $$
  or
  \begin{eqnarray*}&&
  C(n-1)C(n+1)-C^2 (n) > V \\&&
  > C^2 (n)(2d(n)+d^2(n))+C(n-1)C(n+1)(d(n-1)+d(n+1)+d(n-1)d(n+1))\\
  &&>C^2(n)(2d(n)+d^2(n)+d(n-1)+d(n+1)+d(n-1)d(n+1)).
  \end{eqnarray*}
  To satisfy last inequality is sufficient to impose the inequality
  \begin{equation}
  \label{r4}
  C(n-1)C(n+1)-C^2(n) >V >C^2 (n)(2|d(n)| +d^2(n) +|d(n-1)d(n+1)| +|d(n-1)|+|d(n+1)|) 
  \end{equation}
  or
  $$
  \frac{V}{C^2(n)} =\frac{1}{r^2} > 2|d(n)| +d^2(n) +|d(n-1)d(n+1)| +|d(n-1)|+|d(n+1)|.
  $$
  We can assume that $d(m) < 12e^{-r},\ m=n,n+1,n-1.$ Then 
  $$
  2|d(n)| +d^2(n) +|d(n-1)d(n+1)| +|d(n-1)|+|d(n+1)| <48e^{-r}+288e^{-2r} <\frac{1}{r^2},\ r>10.
  $$
  It is left to check inequality~(\ref{l2}) for $r\leq 12,\ n\neq 2,4.$ We can do this with the help of software "Mathematica".
\bigskip

\begin{center}
{\bf Acknowledgment}
\end{center}

\bigskip

We would like to express our gratitude to  Unifesp  and  to USP also where he started this work and specially  Prof. K. Yoshiharu.

\end{document}